\documentclass[authoryear, preprint,12pt]{elsarticle}

\usepackage{amssymb}
\usepackage{amsmath}
\usepackage{makecell}
\usepackage{url}
\usepackage{longtable}
\usepackage{graphicx}
\usepackage{tabularx}
\usepackage{subcaption}
\usepackage{rotating}
\usepackage{booktabs}
\graphicspath{ {./Figures/} }

\journal{Computers \& Chemical Engineering}

\newcommand{\be}{\begin{equation}}
\newcommand{\ee}{\end{equation}}
\newcommand{\bea}{\begin{eqnarray}}
\newcommand{\eea}{\end{eqnarray}}

\begin{document}
\begin{frontmatter}

\title{A Multi-Scale Optimization Framework for Grid-Integrated Electrolysis}

\author{Kiernan X. Jennings, Victor M. Zavala, and Styliani Avraamidou}

\affiliation{organization={Department of Chemical and Biological Engineering, University of Wisconsin-Madison},
            addressline={1415 Engineering Dr}, 
            city={Madison},
            postcode={53706}, 
            state={WI},
            country={USA}}

\begin{abstract}
The increasing penetration of wind and solar resources into the power grid motivates the integration of flexible technologies to dynamically shift power loads in response to grid volatility and emergency events. The water electrolyzer presents a synergistic opportunity to provide flexibility through demand response (DR), while simultaneously electrifying hydrogen production; however, highly dynamic operation schedules accelerate device degradation.
This work presents a mixed-integer linear program (MILP) optimization framework to study the multi-scale coupling between short-term operational flexibility provision in electrolysis devices and long-term stack replacement decisions driven by degradation. 
Active day-ahead market (DAM) participation of a 2.2 MW alkaline water electrolyzer over 22 years is solved as a case study.
Our framework reveals that the multi-scale scheduling of DR operation and replacement decisions can extend optimal stack lifetimes by up to 2 years through load-shifting and further reduce lifetime electricity expenses by 33\% relative to inflexible constant operation. 
Furthermore, we quantify key device parameter tradeoffs and next-generation design goals, where our analysis challenges the feasibility of the standard \$1/kg levelized cost of hydrogen (LCOH) production target solely through market arbitrage. 
Ultimately, this framework quantifies the largely unexploited economic value of multi-scale optimization in grid-integrated electrolysis.
\end{abstract}

    






\end{frontmatter}

\section{Introduction}
The push for decarbonization of the energy sector presents major challenges for the power grid due to the increasing penetration of intermittent wind and solar power. These problems are exacerbated by the escalating growth of large inflexible loads (i.e., data centers). 
These trends necessitate the adoption of more flexible technologies capable of adjusting load and generation in real time to stabilize the grid \citep{bennatoFlexibilityEnergySector2025}. 

One method of providing this flexibility is through demand response (DR). By coordinating large consumers to operate in response to grid signals, DR can offer flexibility to the grid \citep{breeDemandSideManagement2019}. These signals come in the form of electricity market prices. Traditionally, power grid operators coordinate supply and demand of energy systems through hierarchical market transactions in the day-ahead market (DAM) and real-time markets (RTM). A manifestation of increasing intermittency in the supply of energy is the increasing volatility of electricity prices in the DAM and RTM. This opens an opportunity to dynamically schedule electrochemical technologies to provide flexibility for the power grid while performing market arbitrage, subsequently reducing their energy costs \citep{cozzolinoReviewElectrolyzerbasedSystems2024, IndustryProvesBe, eichmanNovelElectrolyzerApplications2014}. Specifically, these devices can conduct electricity market arbitrage by ramping up power consumption at low price periods (low grid congestion) and reduce power usage at high price periods (high grid congestion) \citep{stoustrupSmartGridControl2019, mohanpurkarElectrolyzersEnhancingFlexibility2017}.

We primarily focus on electrolysis devices which have the potential of being key chemical manufacturing and power grid sector-coupling assets \citep{baldeaTransformingProcessIndustries2025}. These electrolysis systems can actively participate in a range of electricity markets \citep{nguyenGridconnectedHydrogenProduction2019, mossleEconomicOptimizationDynamic2025, dowlingMultiscaleOptimizationFramework2017} and have shown this simultaneously minimizes internal operational expenses \citep{maExploitingElectricityMarket2023, eichmanNovelElectrolyzerApplications2014}. However, a major concern with dynamic operation of electrolysis devices is the accumulation of {\em degradation} and the subsequent long-term {\em durability} of these systems under dynamic operation\citep{nezhadkhatamiTechnoeconomicOptimizationPower2026}. 

Short- and medium-term experimental studies (hours to months) on electrolysis devices have shown that irregular and fluctuating operation of electrolysis devices adversely affects durability compared to constant operation \citep{aliaElectrolyzerDurabilityLow2019, jungEffectRampingRate2024}. Specifically, frequent startup and shutdown cycles can shorten the lifetime of the electrolysis stack due to the degradation of cell components \cite{mehdi2023}. Despite this, there are limited studies that analyze the effect of dynamic operation over long time horizons (years to decades) due to the difficulty in conducting experiments \citep{remmeGlobalHydrogenReview2024, weissImpactIntermittentOperation2019, lystbaekReviewEnergyPortfolio2023}. Thus, we must use electrolysis models to be able to capture degradation effects from dynamic operation tractably and accurately. The effects on durability can then inform infrastructure decision making, such as optimal replacement intervals; policy decisions, such as optimal market participation rules; and environmental impacts in life-cycle analyses -- all of which would be poorly characterized if degradation effects are excluded or approximated. Particularly, there has been little to no study on the optimal electrolysis stack replacement resulting from dynamic operation.

The existing models of grid-connected electrolysis devices capture degradation phenomena primarily via heuristics that aim to approximate the accumulation of degradation irrespective of the operational modes. For instance, \citet{parkImpactDegradationEconomics2025} used fixed degradation rates and stack replacement heuristics over a long time horizon. \citet{wangOptimalSchedulingElectrichydrogen2025} quantified degradation as a growing operational expense over time, bypassing the logic of degradation accumulation. \citet{aoOperationStrategyOptimization2025} formulated degradation of the electrolyzer using a fixed linear aging model across the time horizon. 
Given that degradation accumulation is a function of the system's usage, more detailed models are necessary to effectively capture the value of the flexibility provided by electrolysis systems \citep{helistoImpactOperationalDetails2021}.

Another approach to capturing degradation is to model the system's voltage increase as a function of the operation of the electrolysis systems \citep{frenschInfluenceOperationMode2019}, denoted as usage-based degradation. This voltage increase encapsulates thermodynamic penalties from various degradation mechanisms to maintain the same hydrogen production levels. Usage-based degradation incorporates logic to tie together operational protocols and the physical effects impacting durability.  This type of model has been used extensively in applications for optimizing short-term controls and process design. \citet{yangDegradationawareSchedulingThermal2026} constructed a control framework while explicitly accounting for different degradation mechanisms as functions of the operational states. \citet{schofieldDynamicOptimizationProton2024} modeled usage-based degradation for cost optimal design and operation over representative days. Similarly, \citet{filipeMultiObjectiveOptimizationGreen2025} formulated a model for bidding into multiple markets while accounting for usage-based degradation. Conducting years to decades long studies using usage-based degradation functions can be computationally challenging due to the inherent nonlinearity of the electrolysis physics. This limits the ability to study the multi-faceted effects of dynamic operation and the associated degradation.

In this work, we present a multi-scale optimization framework for grid-integrated electrolysis systems, cast as a mixed integer linear program (MILP). Our framework explicitly couples short-term DR operation scheduling under time-varying electricity market prices with usage-based degradation considerations and long-term stack replacement decisions. To the best of our knowledge, this is the first study to co-optimize the full-space hourly operation scheduling with embedded usage-based degradation inside a cohesive multi-year stack replacement decision-making framework. Through this work, we demonstrate that strategically managed DR device operation can {\em reduce cumulative degradation} and {\em extend device lifetimes} under volatile market conditions, contrary to conventional expectations. Ultimately, this framework provides a generalized approach to maximize economic performance of the grid-integrated electrolysis system. This work can be easily extended for projecting economic outcomes of potential electrolysis systems and informing policy decisions. 

The paper is organized as follows. Section \ref{sec:method} introduces the proposed multi-scale optimization framework. In Section \ref{ssec: discretization} and \ref{ssec: super} we discuss the time discretization and the scope of the problem. The optimization model is specified in Section \ref{ssec: model} thoroughly. The novel degradation model and replacement decisions are discussed in \ref{ssec: usage-based} and \ref{ssec: replacement}, respectively. We discuss the case study in Section \ref{sec: case study} and discuss the results in Section \ref{sec: discussion}. We analyze the case study results results in Section \ref{ssec: market} then perform an economic sensitivity analysis in Section \ref{ssec: sensitivity}. Finally, we wrap up with conclusions and future directions in Section \ref{sec: conclusion}. 

\section{Multi-Scale Optimization Formulation}\label{sec:method}

We formulate a MILP model that captures the multi-scale coupling between operational and strategic decision-making for a grid-connected water electrolyzer. We assume that the water electrolyzer tracks electricity market signals of wholesale markets by adjusting its electricity load/demand. We assume that the electrolyzer is a price-taker and produces hydrogen at a fixed selling price.

\subsection{Time Representation}\label{ssec: discretization}

\begin{figure}[!htp]
	\centering
	\includegraphics[width=.7\textwidth]{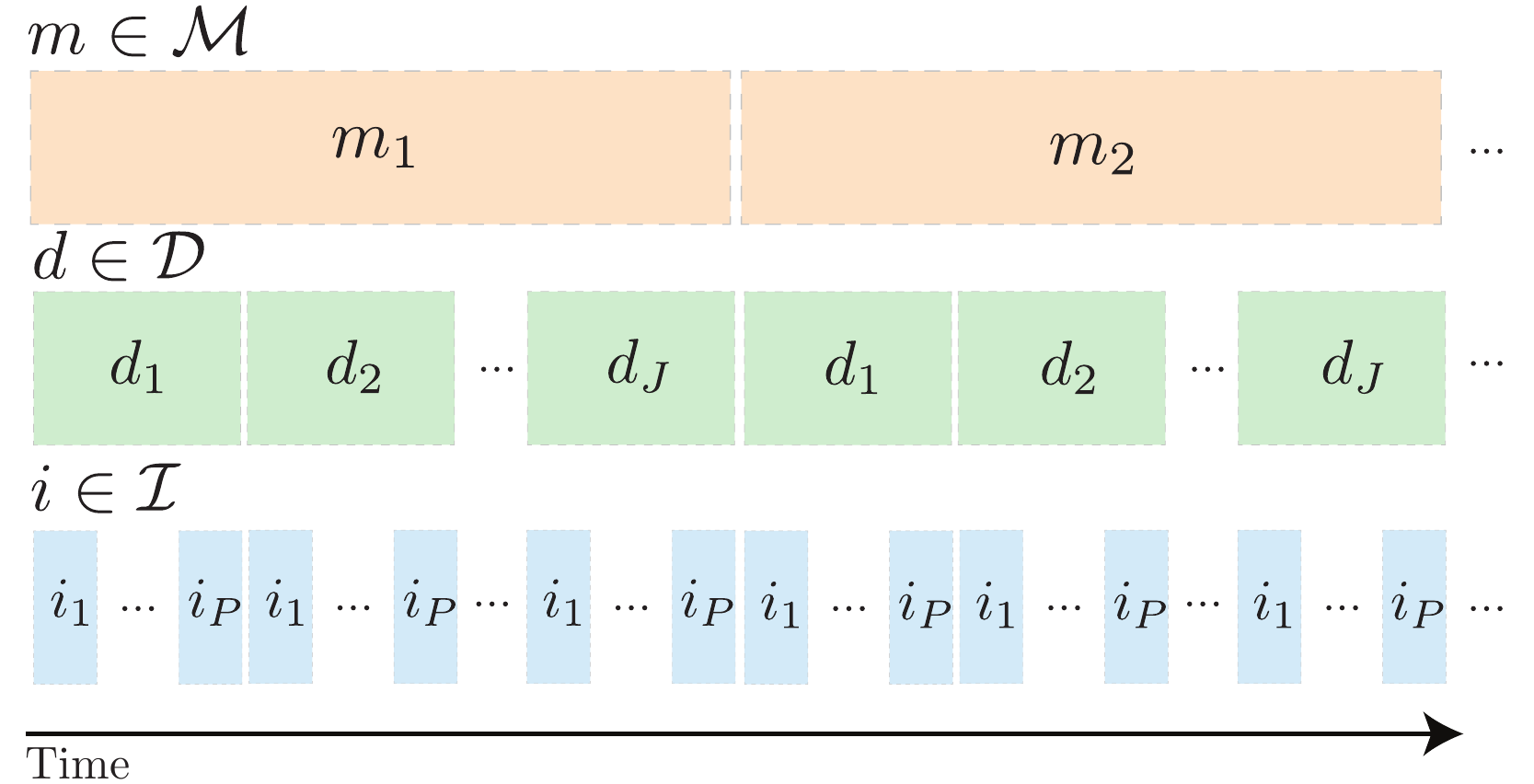}
	\caption{Schematic representation of the time horizons in the multi-scale optimization formulation. For $N$ years, each year ($m$) contains $J$ days ($d$) and each day contains $P$ time instants ($i$).}
	\label{fig: time horizons}
\end{figure}

The time discretization scheme for the model captures the decision making process across multiple, hierarchical timescales (see Figure \ref{fig: time horizons}). Here, we consider annual replacement decisions $(m \in \mathcal{M})$, maintaining a daily chemical demand $(d \in \mathcal{D})$, and market participation at a fixed resolution $(i \in \mathcal{I})$. We assume hourly market participation in the DAM, although the framework can be adjusted to capture participation in RTM. 

We discretize the time horizon of the optimization problem by denoting the lexicographic time set $T$ that captures the hierarchical and multi-scale problem in nature. 
\be
\begin{split}
T &:= \mathcal{M} \times \mathcal{D} \times \mathcal{I} \\
& = \{ (m_1, d_1, i_1), (m_1, d_1, i_2),\dots, (m_{N-1}, d_J, i_P),\\
& (m_N, d_1, i_1), \dots, (m_N, d_2, i_1), \dots,  \\
& (m_N, d_J, i_{P-2}), (m_N, d_J, i_{P-1}), (m_N, d_J, i_P) \}\\
\end{split}
\ee
Each point $(t \in T)$ is defined by a tuple corresponding to a specific time instant and the entire set consists of all time points across the time horizon in chronological order. The notation of $t\text{+}1$ represents a forward time step of length $\gamma$ (e.g., one hour in DAM). Similarly, $t\text{-}1$ represents a backwards time step. For clarity of variables defined on the singular annual timescales, we define a subset of the time points, $T^* \subseteq T$, which contains the first time instant of each year $m$ in the representative long-term timescale. 
\be
T^* := \{(m_1, d_1, i_1), \dots , (m_{N-1}, d_1, i_1), (m_N, d_1, i_1)\}
\ee

\subsection{Model Structure}\label{ssec: super}
\begin{figure}[!htp]
	\centering
	\includegraphics[width=\textwidth]{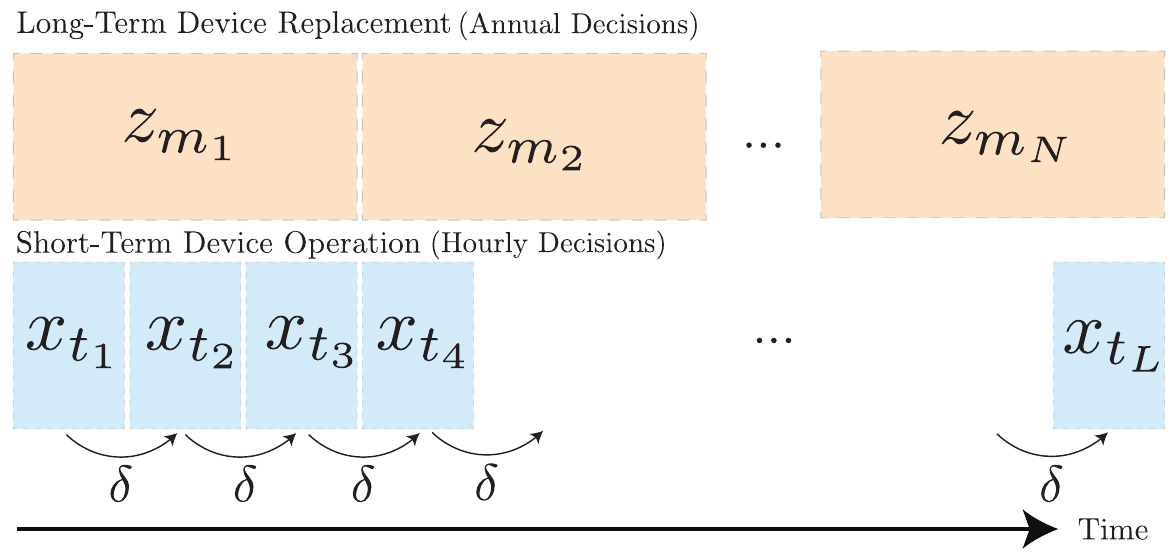}
	\caption{Diagram of the multi-scale optimization structure for decision variables: device replacement ($z_m$) and operation scheduling ($x_t$). Degradation, $\delta$, is accumulated in the short-term operational decision making layer.}
	\label{fig: multiscale diagram}
\end{figure}

The multi-scale structure of the optimization model is represented in Figure \ref{fig: multiscale diagram}, highlighting the coupling of long-term (annual) and short-term (hourly) decisions. The electrolyzer stack replacement decision ($z_m$) is a binary variable that captures replacement of the electrolyzer stack and are made annually. The short-term operational time horizon is composed of binary operation decisions ($x_t$). These decisions are for each operation mode 'on', 'off', and 'standby' and made for each hour. More details on this are found in \ref{sec: operation modes}.

\subsection{Optimization Model}\label{ssec: model}

\subsubsection{Economic Objectives}
The objective function is defined as the maximization of the net present value (NPV) of the electrolyzer over the entire time horizon. 
This objective accounts for stack replacement expenses ($REPEX_m$), short-term operational expenses ($OPEX_m$) and hydrogen production revenue ($REV_m$) for a given discount rate, $\rho$. The initial capital expenditure ($CAPEX$) is assumed to be a fixed expense, dependent on the system of interest.
\be\label{eq: obj fn}
\max_{\textbf{z,x}} NPV = - CAPEX + \sum_{m \in \mathcal{M}} \frac{1}{(1+\rho)^m} \left( REV_m - REPEX_{m} - OPEX_m \right)
\ee

The capital expenditure is set to the parameter, $\psi^{CAPEX}$, and problem specific. We assume this expense to be made in year 0 with no discount rate.
\be
CAPEX = \psi^{CAPEX}
\ee

The stack replacement cost for each year is calculated as:
\be
REPEX_m = \psi^{stack} z_m \quad \forall m \in \mathcal{M}
\ee
for a fixed price of the electrolyzer stack, $\psi^{stack}$, over the entire time horizon.

The operational expenses for each year are calculated based on fixed operational expenses, $\psi^{OPEX}$, and the variable expenses from market participation. We assume the electricity prices, $\psi^E_{m,d,i}$, to be constant for each time $t$.
\be
OPEX_m = \psi^{OPEX} + \sum_{d \in \mathcal{D}}\sum_{i \in \mathcal{I}} \psi^E_{m,d,i} e_{m,d,i} \quad \forall m \in \mathcal{M}\\
\ee
The fixed operational expenses, $\psi^{OPEX}$, are calculated as a fraction of the initial capital expenditure on a yearly basis \citep{namiTechnoeconomicAnalysisCurrent2022}.

The revenue from the electrolyzer comes from selling the hydrogen at a fixed market price, $\psi^{H}$. 
\be
REV_m = \sum_{d \in \mathcal{D}}\sum_{i \in \mathcal{I}} \psi^H h_{m,d,i} \quad \forall m \in \mathcal{M}
\ee
Adaptations of this formulation can include additional forms of electricity market remuneration (e.g., regulation payments).

The levelized cost of hydrogen (LCOH) is used to compare economic feasibility of this operation. This value is not optimized over but rather provides a basis to compare different hydrogen production technologies.
\be 
LCOH = \frac{\psi^{CAPEX} + \sum_{m \in \mathcal{M}}(1+\rho)^{-m}(REPEX_{m} + OPEX_m)}{\sum_{m \in \mathcal{M}}(1+\rho)^{-m} \sum_{d \in \mathcal{D}} \sum_{i \in \mathcal{I}} h_{m,d,i}}
\ee
This calculation provides the cost of producing hydrogen in USD per kilogram.

\subsubsection{Operation Modes}\label{sec: operation modes}

Operation decisions are modeled as binary variables $x_t$ for each operation mode at time $t$. 
The electrolyzer operational model by \citet{baumhofOptimizationHybridPower2023} is used in this formulation. The electrolyzer is modeled with three operational states, 'on', 'off' and 'standby'. This choice is modeled as using the following constraint:
\be
x^{on}_{t} +	x^{off}_{t} +	x^{sb}_{t} = 1 \quad \forall t \in T
\ee 
The 'on' state represents operating the electrolyzer at the nominal power rating to produce hydrogen. 
The 'off' state represents a cold shut off of the electrolyzer system with no power consumption or chemical production. 
The 'standby' mode offers the decision to operate at a fractional amount of energy to prevent thermal losses in the electrolyzer system, but produce no hydrogen.

An additional binary variable, $x^{start}_t$, is introduced to indicate whether there is a transition from the operational 'off' mode to 'on' mode occurring at time $(t-1)$ to time $t$. This will be used to track the number of cold starts in the system and incorporate the associated degradation to the stack \cite{mehdi2023}.
\begin{subequations}
	\begin{align}
		&x^{start}_{t} \geq x^{on}_{t} - x^{on}_{t-1} - x^{sb}_{t-1} &  \forall t \in T\setminus \{t_1\} \\
		&x^{start}_{t} \leq x^{on}_{t} & \forall t \in T\setminus \{t_1\} \\
		&x^{start}_{t} \leq 1 - x^{on}_{t-1} - x^{sb}_{t-1} & \forall t \in T\setminus \{t_1\}
	\end{align}
\end{subequations}
We do not track the other transition of 'standby' to 'on' as we assume that the warm standby mode prevents the degradation from cold starts of the electrolyzer system. 

\subsubsection{Physical Model}\label{sec: degrad}

The electrolysis physics underlying the operation of an electrolysis system is nonlinear in nature which would significantly impact the tractability of this multi-scale framework. Thus, a series of linearizations and reformulations are used to formulate the original mixed-integer nonlinear program into a MILP.

A common linearization is the hydrogen production curve for operation at the nominal power rating of the electrolyzer \cite{varelaModelingAlkalineWater2021}. The slope of this hydrogen production linear approximation is the efficiency of the electrolyzer $a_t$ at time $t$. This efficiency is given at a rate of hydrogen mass produced per power consumed. This rating is analogous to the percentage of the lower heating value (LHV) of hydrogen.  The hydrogen production at time $t$ is dependent on the operation decision of the electrolyzer and the present efficiency of the device, $a_t$. 
\be 
h_{t} = \phi \gamma \  a_{t} x_{t}^{on} + \beta \gamma \ x_{t}^{on} \quad \forall t \in T.
\ee
The hydrogen production is also a function of the nominal capacity of the electrolyzer stack, $\phi$, the time discretization of the model, $\gamma$.
The intercept of the hydrogen production curve, $\beta$, is treated as constant across the entire time horizon in units of kg H$_2$/hr.

In this implementation, the production expression is nonlinear with respect to the $a_{t}$ and $x_t$ terms. Thus, an exact linearization term is used: 
\be
h_{t} = \phi \gamma \  w_{t} + \beta \gamma \ x_{t}^{on} \quad \forall t  \in T,
\label{eq: hydrogen production}
\ee
where, 
\begin{subequations}
	\begin{align}
		&w_t \leq \bar{\alpha}  x^{on}_t & \forall t \in T,\\
		&w_t \geq \underline{\alpha}  x^{on}_t & \forall t \in T,\\
		&w_t \geq a_t - \bar{\alpha} (1 - x^{on}_t) & \forall t \in T,\\
		&w_t \leq a_t - \underline{\alpha} (1 - x^{on}_t) & \forall t \in T,
	\end{align}
\end{subequations}
The upper and lower bounds on the efficiency state, $\bar{\alpha}, \underline{\alpha}$, are used in the exact McCormick envelope. A general guideline and heuristic commonly used is a stack lifetime of 10 years, which equivalently can serve as a generic lower bound \citep{TechnicalTargetsLiquid2025}.

The power demand of the electrolyzer is defined based on the fixed capacity of the system and the operation decisions. Given the three state system defined, the total amount of power comes from either the on or standby operation modes.
\be
e_{t} = \phi \gamma \ x^{on}_{t} + \phi \gamma \zeta^{sb} \ x^{sb}_t  \quad \forall t \in T.
\ee
where $\zeta^{sb}$ is the percentage of the nominal load to run in standby mode, typically around 1-5\% \citep{varelaModelingAlkalineWater2021}.

Similarly for finer time discretizations, $\gamma$, a ramping limitation, $\theta$, of the device is necessary to describe the inability to transition to any power load instantaneously. This varies based on the type of electrolyzer technology being modeled. 
\be
e_t - e_{t-1} \leq \theta \phi \quad \forall t \in T \setminus \{t_1\}
\ee
A typical assumption excludes this constraint as the ramping rates for electrolyzers lie around $\pm20\%$ nominal power per second \citep{baumhofOptimizationHybridPower2023}. Given hourly operation time discretization when participating in the DAM, ramping limitations can be neglected for the device (e.g., $\theta=1$) \citep{varelaModelingAlkalineWater2021}.

\subsubsection{Chemical Demand}

We enforce meeting a minimum operation quota and subsequent hydrogen demand of production from the electrolyzer. The medium-term timescale $d \in \mathcal{D}$ groups together time instants over the course of the day. We then enforce meeting the minimum daily chemical demand, $\sigma$, in terms of mass of hydrogen produced.  This can easily be extended to be defined on weekly or monthly demands.
\be
\sum_{i \in \mathcal{I}} h_{m, d, i} \geq \sigma \quad \forall d \in \mathcal{D}
\ee

\subsection{Usage-Based Degradation Model} \label{ssec: usage-based}

The degradation of the electrolyzer device is treated as a fixed efficiency loss on the electrolyzer stack as a function of the binary operational decisions. The degradation state acts as a coupling constraint across the entirety of the time horizon. We consider two degradation types: operational and startup. Operational degradation relates to the decline in electrolyzer efficiency due to continuous operation conditions on a per hour basis. The startup degradation comes from the 'off' to 'on' operation of electrolyzers and represents a cold start. 
\be 
	a_{t} = a_{t-1} - \gamma \delta^{op} x_{t}^{on} - \delta^{start} x_{t}^{start} \quad \forall t \in T\setminus (\{t_1\} \cup T^*)
\ee 
This distinction in operational and startup degradation types allows us to draw comparisons across constant operation against dynamic operation. We note that the startup degradation is avoided when going from 'standby' to 'on' operation modes. A third degradation mode from 'standby' to 'on' is neglected for tractability as it would remain significantly smaller than the cold startup degradation \citep{martinezlopezDynamicOperationWater2023}.

Degradation for electrolysis is commonly experimentally characterized by the increase in the stack voltages. This increase in voltage, $\eta$, can lead to an equivalent loss of production efficiency for the device.
\be\label{eq: degradation}
\delta = \frac{V_0}{\text{HHV} } \frac{\eta}{V_i (V_i + \eta)}
\ee
This equates the hydrogen production efficiency loss for a given voltage degradation increase of both degradation types, operational and startup. The thermoneutral voltage of water electrolysis, $V_0$, represents the ideal operation voltage for isothermal operation. The higher heating value of hydrogen, $HHV$, is used to equate the increase in voltages to the production rate of hydrogen from the system. The initial voltage of the electrolysis cell, $V_i$, encapsulates all types of activation, ohmic, and transport overpotentials. This formulation for the efficiency decrease parameter is used for the entirety of the optimization problem, thus the associated degradation value for operation and startup is assumed to be a constant parameter over the entire of the time horizon.

\subsection{Stack Replacement Decisions}\label{ssec: replacement}

Optimal replacement of the electrolyzer stack is used to bridge the gap between operational modes and accumulative degradation and the long-term infrastructure decisions. To model replacement, the electrolyzer efficiency is formulated to reset back to the initial efficiency of the electrolyzer stack, $\alpha^0$. The binary variable, $z$, denotes decision to replace in a given year $m$. 
We formulate this at the specified intervals of $t \in T^*$ at the initial time instant of each year.
\be 
\begin{split}
	a_t &= a_{t-1} - \gamma \delta^{op} x_{t}^{on} - \delta^{start} x_{t}^{start}  +(\alpha^0 - a_{t-1} )z_{m} \quad \forall t \in T^*, m \in \mathcal{M}
\end{split}
\ee

This constraint is nonlinear with respect to the $a_{t-1}$ and $z_m$ terms. The same exact McCormick linearization from the hydrogen production is used.
\be
a_{t} = a_{t-1} - \gamma \delta^{op} x_{t}^{on} - \delta^{start} x_{t}^{start} + \alpha^0 z_{m} - v_{m} \quad \forall t \in T^*, m \in \mathcal{M}
\ee
where we introduce an auxiliary variable $v_m$ which is constrained by the following:
\begin{subequations} 
	\begin{align}
		&v_m \leq \alpha^0 * z_m& \forall m \in \mathcal{M}  \\
		&v_m \geq \underline{\alpha} * z_m &\forall m \in \mathcal{M}  \\
		&v_m \geq a_{t-1} - \alpha^0 (1 - z_{m}) & \forall t \in T^*, m \in \mathcal{M}  \\
		&v_m \leq a_{t-1} - \underline{\alpha}(1 - z_{m}) & \forall t \in T^*, m \in \mathcal{M} 
	\end{align}
\end{subequations}
Here, we assume electrolyzer stack replacements occur in one time instant with no device downtime.

\section{Case Study}\label{sec: case study}

To illustrate the usability of the proposed formulation, we consider an alkaline water electrolyzer (AWE) participating in the Electric Reliability Council of Texas (ERCOT) independent system operator (ISO) region in the DAM. The AWE is a mature electrolysis technology that has been commercialized and implemented to large scale in capacities exceeding 20 MW. This type of electrolysis system includes coupled electrodes operating in an alkaline electrolyte separated by a diaphragm to produce hydrogen and oxygen gas. The benefit of this type of technology is the lower and more predictable rate of degradation despite a slower ramping capability \citep{campbell-stanwayTechnoeconomicAnalysisElectrolyser2025}.

We use industrially relevant AWE device parameters (see in Appendix Table \ref{tbl: device params}) to analyze the long-term performance and durability when actively participating under DR.
Texas was selected as a representative location for the case study based on high regional demand of hydrogen production for the chemical industry \citep{hickeyLeveragingGulfCoast2024}. The historical electricity price data is sourced from the ERCOT DAM and the annual Texas average industrial distribution network electricity prices \citep{electricreliabilitycounciloftexasercotDAMHourlyLMPs2024, mcgrathElectricPowerMonthly2026}. 
Electricity market price datasets for 2014-2024 are replicated to generate a 22-year time horizon at an hourly resolution to demonstrate the capabilities of the modeling framework over long time horizons. Future work can extend to using price forecasting to project grid-integrated electrolysis performance.
The startup voltage degradation was determined from experimental cycling protocols in \citet{parkImpactDegradationEconomics2025} and calculated from the average voltage increase across 500 cycles. The operational voltage degradation was from U.S. Department of Energy 2022 standards for AWE and consistent with reported values in literature \citep{TechnicalTargetsLiquid2025}.

\begin{figure}[!htp]
    \centering
    \includegraphics[width=0.8\textwidth]{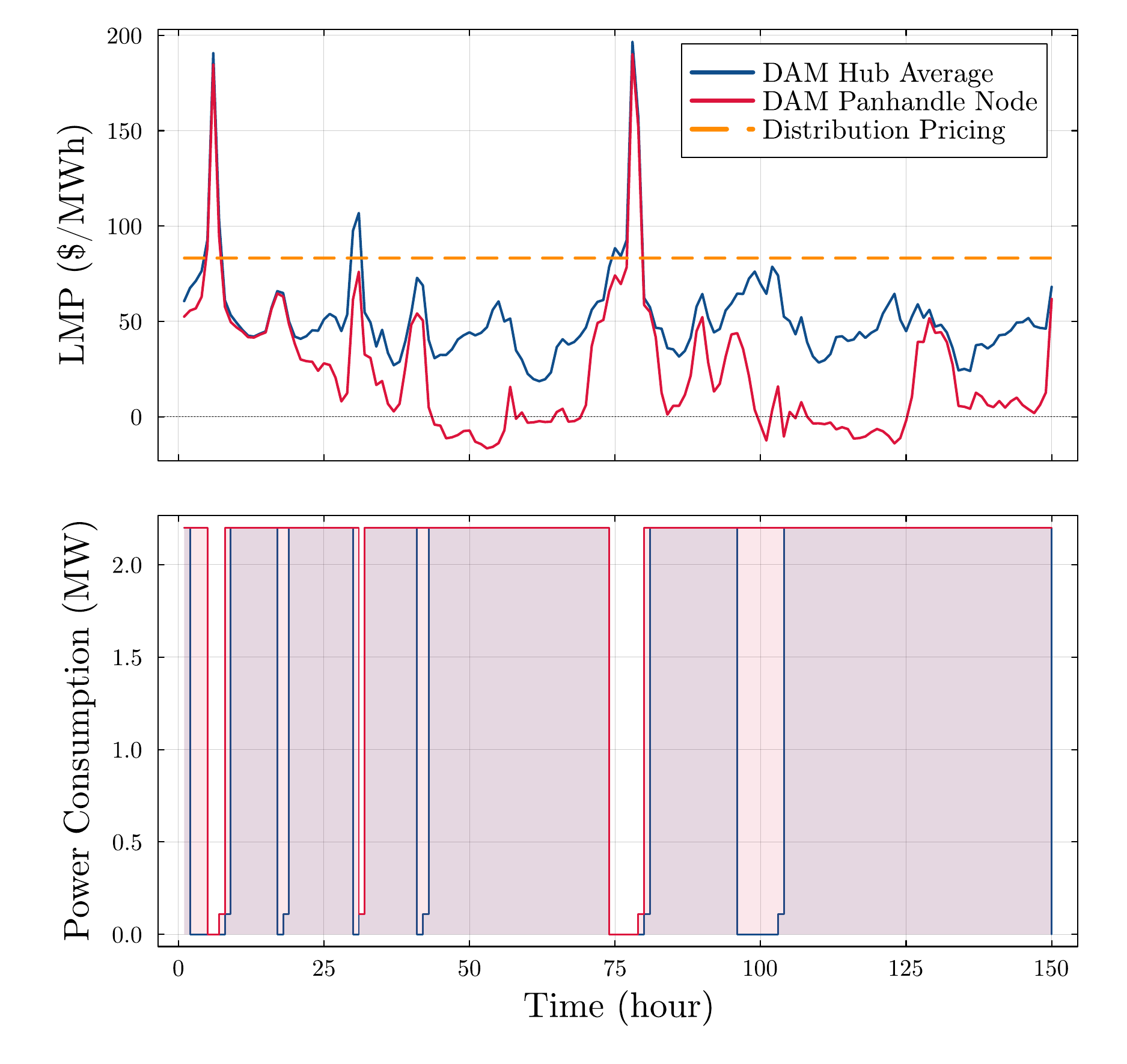}
    \caption{Top: Locational marginal price (LMP) from the DAM hub average, panhandle node, and distribution pricing for a representative week in April 2022.
    Bottom: Power consumption for the flexible operation of 2.2 MW AWE for given price profile from above.}
    \label{fig: price profiles}
    \label{fig: flexible operation}
\end{figure}

We compare three different price time profiles for market participation in Figure \ref{fig: price profiles} and Table \ref{tbl: ercot characteristics}. The first profile is the hub average locational marginal pricing (LMP) in the ERCOT ISO. This represents the average market price across all nodes in the ERCOT domain. The second is the ERCOT panhandle LMP profile, which is located adjacent to high wind power generation (Post Wind Farm LP). These two profiles were selected to assess the potential benefits of participation adjacent to renewable power resources. The third profile is the average annual distribution electricity price. This is treated as constant over the calendar year. Participation in this regime tends towards market agreements between the consumer and the municipal electricity providers.

\begin{table}[!htp]
	\centering
	\caption{Comparison of electricity markets statistics for 2014-2024. The DAM pricing profiles were taken as time series data from the ERCOT ISO. The distribution pricing profile for the Texas case study was taken as a constant across the entirety of each year. High prices denoted as $>$ \$250/MWh and extreme prices denoted as $>$\$1000/MWh.}
	\label{tbl: ercot characteristics}
	\small
	\begin{tabularx}{\textwidth}{X c c c}
		\toprule
		& \textbf{Hub Average} & \textbf{Panhandle Node} & \textbf{Distribution} \\
        \textit{Price Statistics} & (\$/MWh) & (\$/MWh) & (\$/MWh) \\
		\midrule
		Average Price & 45.20 & 42.90 & 72.30\\
		Median Price & 23.50 & 21.90 & 69.10\\
		Std. Dev. Price & 286.90 & 287.50 & 5.90\\
		Minimum Price & -2.30 & -36.60 & 66.70\\
		Maximum Price & 8,995.70 & 9,001.00 & 83.20\\
		\midrule
		\textit{Outlier Metrics (count)} \\
		Negative Prices  & 11 & 2,188 & 0 \\
        High Prices  & 862& 872& 0 \\
        Extreme Prices & 308 & 310 & 0\\
		\bottomrule
	\end{tabularx}
\end{table}

Additionally, we differentiate two types of electricity market participation strategies.
Flexible operation denotes DR operation from the presented optimization framework and non-flexible operation denotes constant operation irrespective of the electricity market prices. The latter is representative of traditional chemical manufacturing procedures. An example excerpt of flexible operation is seen in Figure \ref{fig: flexible operation} showcasing active electricity market participation in the hub average price profile. The power consumption of the electrolyzer device fluctuates based on the grid signals (LMP) at a given location. Here, the device selectively decides when to produce at the nominal power capacity and shutdown during high price periods. 

\section{Results and Discussion}\label{sec: discussion}

The resulting MILP for this case study is solved using JuMP modeling language \citep{Lubin2023} with the Gurobi 13.0 solver \citep{gurobi}. The 22-year time horizon model contains 1,539,132 variables (769,566 binary) and 2,124,345 linear constraints. All results are generated on a MacBook Pro M4 with 36 GB of RAM to a $<$1\% MIP gap. 

\subsection{Market Participation} \label{ssec: market}

\begin{figure}[!htp]
	\centering
	\includegraphics[width=\textwidth]{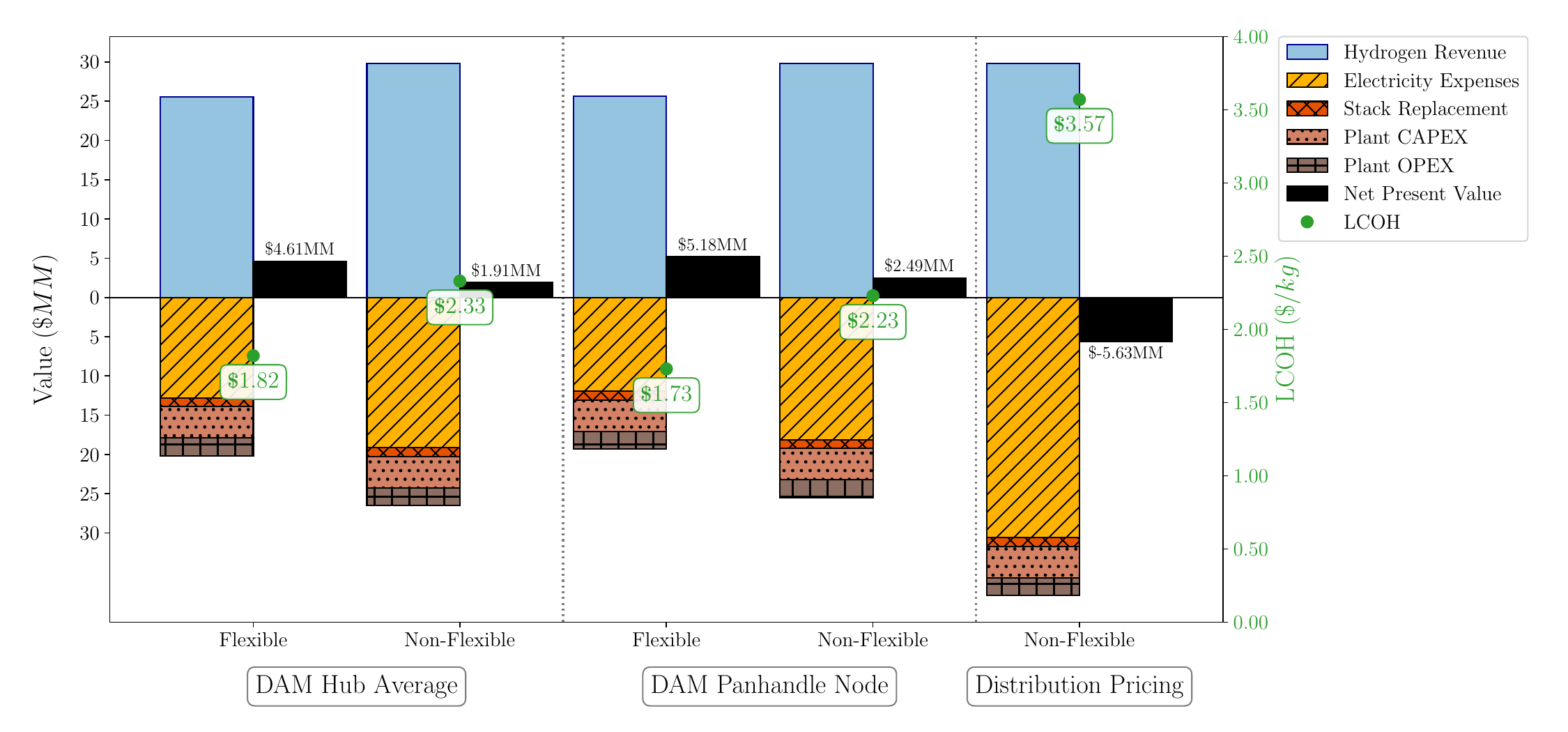}
	\caption{2.2 MW AWE performance in flexible and non-flexible market participation scenarios. Three markets include ERCOT hub average prices, ERCOT panhandle node prices, and average industrial distribution electricity prices over 22-year time horizon. Positive and negative value metrics represent revenue and expenses, respectively.}
	\label{fig: ercot_results}
	\label{fig: electricity expenses}
\end{figure}

\begin{table}[!htp]
	\centering
	\caption{Comparison of participation scenarios statistics for the AWE in Texas. The percentage of operation modes derived as percentage of hours across the entire 22-year time horizon (192,386 hours total). The number of starts and replacement are similarly defined across the time horizon. Non-flexible operation statistics remain the same across all market price profiles.}
	\label{tbl: result summary}
	\small
	\begin{tabularx}{\textwidth}{X c c c} 
    \toprule
		& \multicolumn{2}{c}{\textbf{Flexible}} & \textbf{Non-Flexible} \\
		\cmidrule(r){2-3} \cmidrule(l){4-4}
		\textbf{Price Profile} & \textit{Hub Average} & \textit{Panhandle Node} & \textit{All}\\
		\midrule
		\multicolumn{4}{l}{\textit{Operation Mode (hours, \% total)}} \\
	    On Operation & 182,296 (94.76\%) & 183,171 (95.21\%) & 192,386 (100\%) \\
		Off Operation & 7,529 (3.91\%) & 6,841 (3.56\%) & 0 (0\%)\\
		Standby Operation  & 2,561 (1.33\%) & 2,374 (1.23\%) & 0 (0\%)   \\
		\midrule
        \multicolumn{4}{l}{\textit{Production Statistics}} \\
		H$_2$ Production (t)  & 9,134.2 & 9,176.4 & 9,651.4 \\
        Electricity Consumption (GWh) & 401.3 & 403.2 & 423.2 \\
        \midrule
		\multicolumn{4}{l}{\textit{Device Wear (count)}} \\
		Cold Starts & 161 & 102 & 0 \\
		Warm Starts & 1,879 & 1,906 & 0 \\
		Stack Replacements & 2 & 2 & 2 \\
		\bottomrule
	\end{tabularx}
\end{table}

The results across the five representative participation scenarios in the case study are shown in Figure \ref{fig: ercot_results} and Table \ref{tbl: result summary}. A flexible distribution pricing case was excluded as arbitrage under static pricing offers no benefit. Generally, flexible market participation performed the best, providing at least \$2.7M higher net present value across the 22-year time horizon compared to the respective non-flexible market participation. Industrial distribution pricing was the costliest method of purchasing power, leading to an overall negative NPV of -\$5.63M. Besides this, the cheapest method of purchasing power came from participation in the ERCOT ISO at the panhandle node with higher overall NPV of at least \$0.58M compared to the hub average counterparts. Given the assumption of perfect foresight, this NPV represents an ideal upper bound on device performance.

For the flexible DR operation mode, the device optimally balances the current cost of electricity and the respective efficiency of the device for a given time with foresight of the future market. Flexible market participation contributed a \$0.50/kg reduction in the LCOH compared to non-flexible static participation. Remarkably, the electrolyzer operates about 95\% of the time for both ERCOT LMP profiles analyzed. Despite this, the electrolyzer is able to reduce the electricity expenses by over 33\% compared to the non-flexible operation (see Figure \ref{fig: electricity expenses}).
Irrespective of the electricity market characteristics, flexibility provision fundamentally reduced operation expenses. This subsequently increases hope for industrial viability of DR operation schemes.

A notable difference between the flexible and non-flexible operation modes was the total amount of hydrogen produced over the time horizon. 
The flexible operation produced around 14\% less hydrogen compared to non-flexible participation despite maintaining operation over 95\% of the time horizon. This indicates that the trade-off for deploying grid-integrated electrolyzer devices would involve over-sizing the device to compensate to the same amount of total hydrogen production for non-flexible devices.

We note a minor difference on the economic performance with respect to the location of the electrolysis system, despite the varying levels of renewable intermittent energy penetration. The profitability of the electrolyzer placed adjacent to wind farms only contributed to a total NPV difference of \$570k over the course of the 22 years. This can be attributed to the lower overall average energy price and the corresponding increase in operation based on this price difference. From Table \ref{tbl: result summary}, the electrolyzer operated 0.45\% longer in the panhandle LMP profile compared to the hub average price profile. These results indicate there is only a minor regional impact on the economic performance of grid-integrated electrolysis systems. Further analysis of long-term degradation in hybrid systems (directly interconnected electrolysis systems to renewable electricity generation) is left for future work.

\begin{figure}[!htp]
	\centering
	\includegraphics[width=0.7\textwidth]{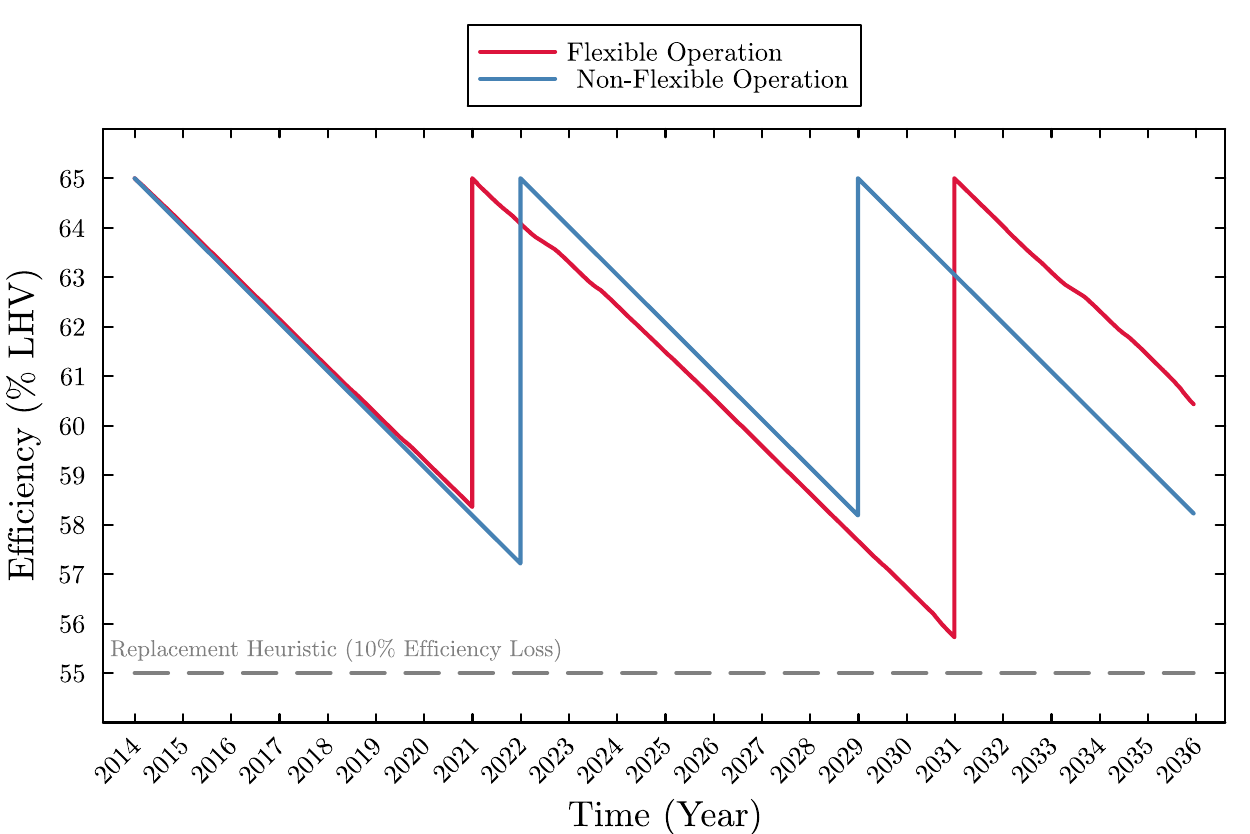}
	\caption{Comparison of the overall efficiency of the electrolyzer over time for flexible and non-flexible participation in the ERCOT hub average profile. Replacements of the electrolyzer are modeled as renewals to the initial efficiency state in a given year. The same trend follows for ERCOT panhandle pricing profiles.}
	\label{fig: replacement}
\end{figure}

\begin{figure}[!htp]
    \centering
    \includegraphics[width=0.7\linewidth]{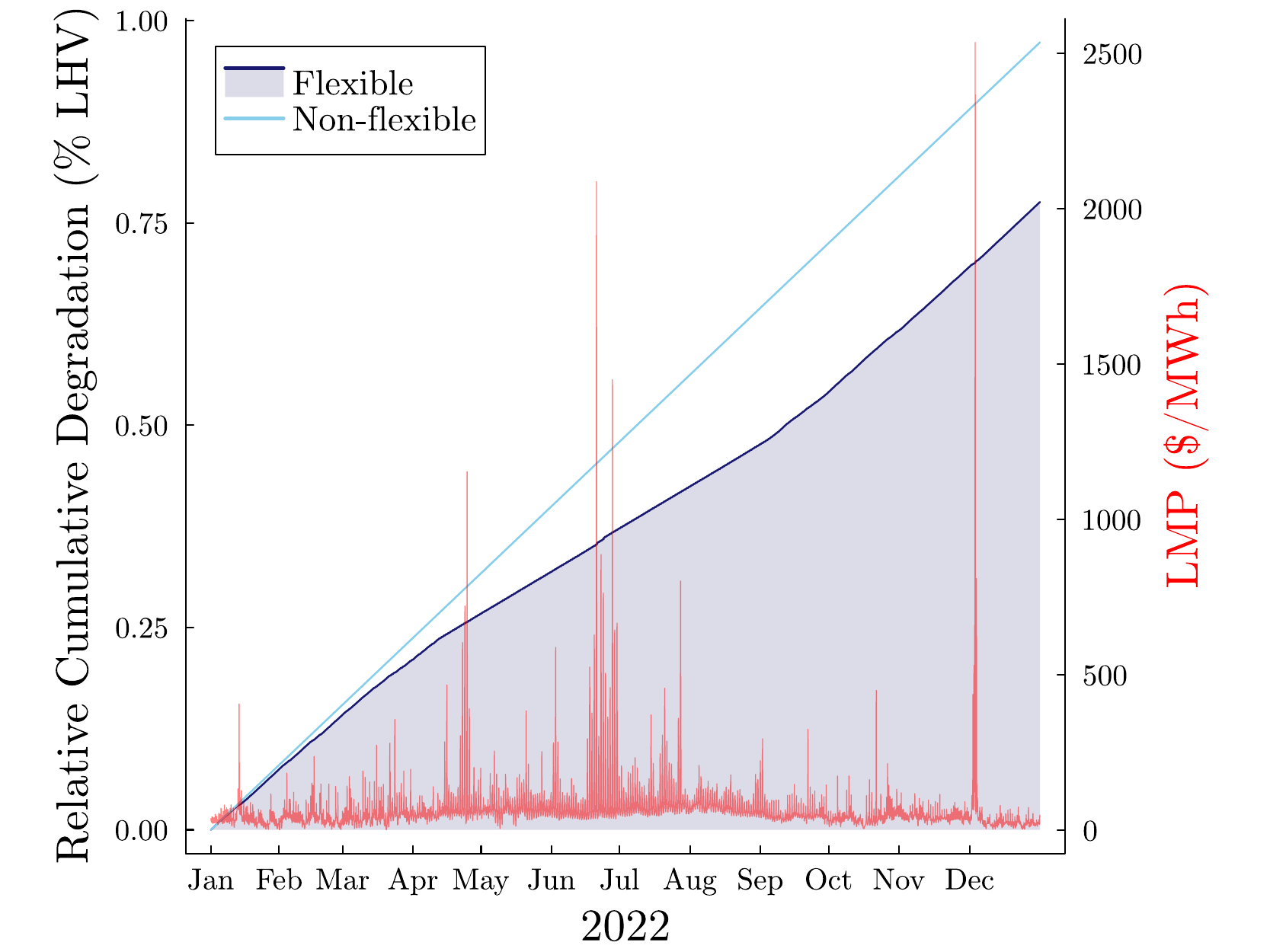}
    \caption{Comparison of cumulative efficiency loss for both flexible and non-flexible operation schemes in the ERCOT hub average profile. Trend follows for ERCOT panhandle pricing profiles. The flexible operation reduces the overall degradation accumulation with intermittently available energy.}
    \label{fig: flex vs nonflex}
\end{figure}

Figure \ref{fig: replacement} compares the rate and frequency of replacement across the flexible and non-flexible operation profiles. Despite formulating the startup degradation mechanism, the frequency of replacements remained constant for all pricing and operation scenarios at two over 22-years. However, the stack lifetimes varied between flexible and non-flexible operation. The former showed stack lifetimes of 7 and 10 years. The latter replaced consistently at stack lifetimes of 7-8 years. These results indicate that the optimal operation and replacement schemes thrive in highly volatile market conditions compared to static, non-flexible operation. Optimal stack lifetimes are prolonged during highly volatile market periods (2021-2025). This phenomenon is also noted in Figure \ref{fig: flex vs nonflex}. These results promise positive capabilities for grid-connected electrolyzer devices as renewable resource penetration continues to increase.

Figure \ref{fig: flex vs nonflex} shows the year 2022 as a representative sample of a highly volatile market and a comparison of the cumulative degradation for flexible and non-flexible operation. Specifically, the 2022 market volatility can be attributed to the tremendous increase (9 TWh, \citep{2022StateMarket2023}) in solar production capacity at the disposal of ERCOT.

Our results demonstrate that as the market prices reach peak volatility in the summer months, flexible operation of the electrolyzer system can reduce the accumulated absolute efficiency degradation by 0.25\% LHV (0.075 kg/MWh) in one year. The electrolyzer optimally sheds operating hours to reduce electricity expenses during peak hours and simultaneously avoids the high rates of degradation faced with cold starts. The effects of this across the given 22-year time horizon enables extending the optimal lifetime of the electrolyzer stack, seen from Figure \ref{fig: replacement}. 

These results illustrate the relevance of capturing the multi-scale coupling of operation and replacement. With previous work (\citep{parkImpactDegradationEconomics2025, aoOperationStrategyOptimization2025}), replacement decision making is fixed over the time horizon and is incapable of dynamically responding to market conditions. In highly volatile markets, we show that DR operation can {\em extend the lifetime} of electrolysis systems when participating optimally with respect to power grid price signals relative to constant operation schemes. Further extensions of this case study may include using scenario-based price forecasts to reinforce our conclusion.

\subsection{Sensitivity Analysis}\label{ssec: sensitivity}
Given the uncertain nature of electrolysis device parameters, we perform a sensitivity analysis to determine the impact of different device and economic parameters on the overall NPV of the electrolyzer across the entire time horizon. First, we seek to determine the effect on perturbations of device parameters to find the optimal parameters to target for electrolyzer device goals. Second, we analyze specified parameters to construct a parameter heatmap design tradeoff between device parameters. 

We run the multi-scale framework for 11 years considering the same 2.2 MW AWE with the baseline parameters (see Appendix Table \ref{tbl: device params}). We chose this planning horizon for the tractability of degenerate problem instances. Extensions for decomposing this framework is left for future work. We consider the flexible DR operation case of the electrolyzer given the hub average LMP pricing profile to be the standard for the sensitivity analysis.

\subsubsection{Design Goals}

\begin{figure}[!htp]
	\centering
	\includegraphics[width=0.7\textwidth]{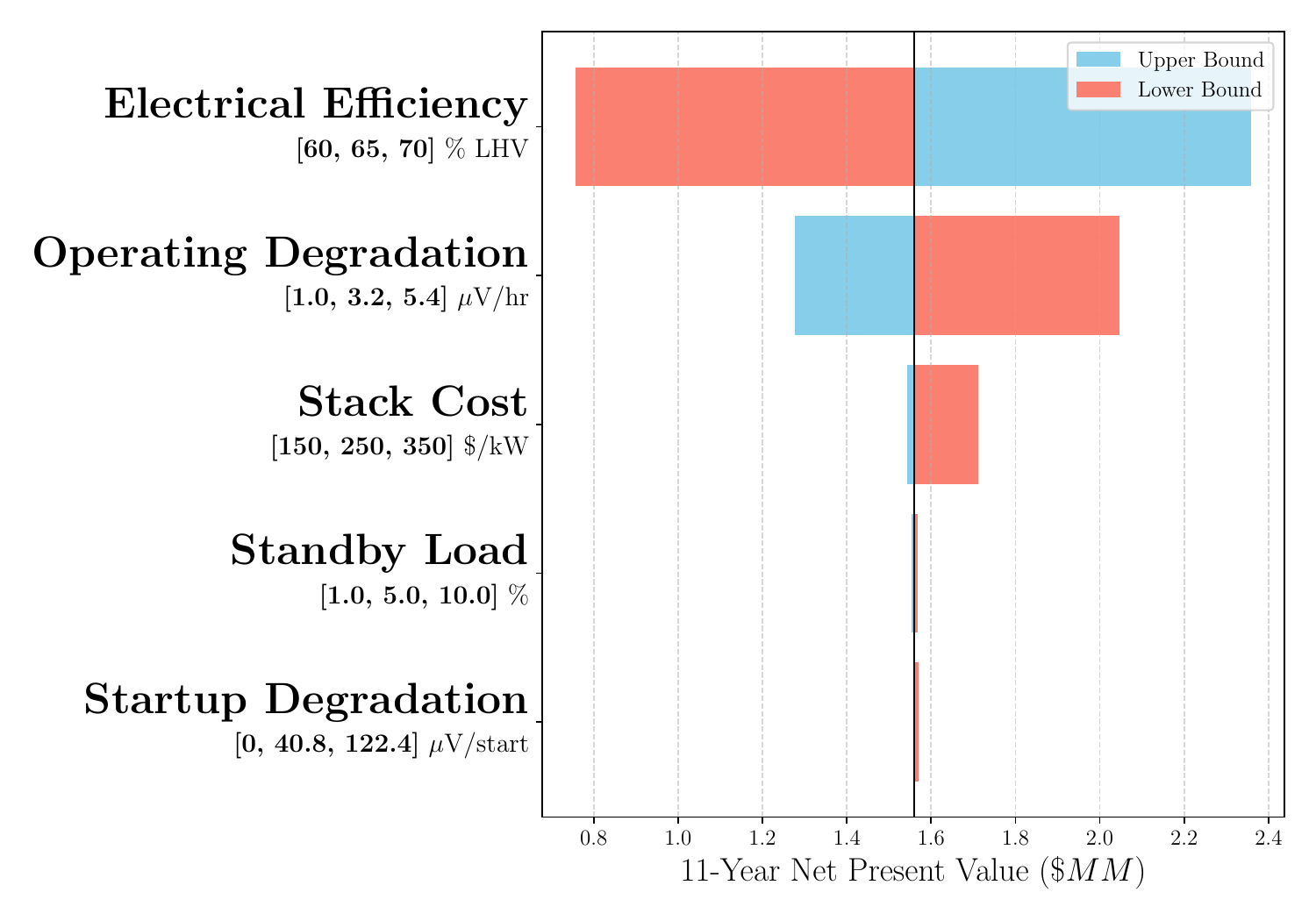} 
	\caption{Impact of electrolysis device and economic parameters on the NPV of the electrolyzer over the 11-year time horizon. Parameter ordering: [Lower Bound, Baseline, Upper Bound]. The baseline NPV is \$1.56MM.}
	\label{fig: tornado plot}
\end{figure}

Figure \ref{fig: tornado plot} compares the economic effects of perturbed electrolyzer device parameters over a 11-year time horizon considering the hub average LMP price profile from Section \ref{sec: case study}. The operation and startup degradation, electrolyzer stack replacement costs, standby percent load, and the initial efficiency of the electrolyzer device are all perturbed to observe the economic value. The range for each parameter comes from technical goals across literature and the DOE \citep{limImpactVoltageDegradation2021, TechnicalTargetsLiquid2025, krishnanPresentFutureCost2023}. 

\begin{table}[!htp]
    \centering
    \caption{Tornado sensitivity analysis results -- Part 1. The three most impactful parameters (largest NPV range) displayed here.}
    \label{tbl: tornado part1}
    \small
    \begin{tabularx}{\textwidth}{X cc cc cc}
        \toprule
        \textbf{Metric} & \multicolumn{2}{c}{\textbf{Electrical}} & \multicolumn{2}{c}{\textbf{Operation}} & \multicolumn{2}{c}{\textbf{Stack Cost}} \\
        & \multicolumn{2}{c}{\textbf{Efficiency}} & \multicolumn{2}{c}{\textbf{Degradation}} & & \\
        \cmidrule(r){2-3} \cmidrule(r){4-5} \cmidrule(r){6-7}
        & \textit{60\%} & \textit{70\%} & \textit{1.0} & \textit{5.4} & \textit{150} & \textit{350} \\
        & \multicolumn{2}{c}{\%LHV} & \multicolumn{2}{c}{$\mu$V/hr} & \multicolumn{2}{c}{\$/kW} \\
        \midrule
        \multicolumn{7}{l}{\textit{Economic Performance}} \\
        NPV (\$MM) & 0.76& 2.36 & 2.05 & 1.28 & 1.71 & 1.54 \\
        LCOH (\$/kg) & 2.15 & 1.86 & 1.82 & 2.04 & 1.95 & 1.92 \\
        \midrule
        \multicolumn{7}{l}{\textit{Operation Mode (hours, \% total)}} \\
        On & 91,094 & 91,852 & 87,232 & 91,115 & 91,500 & 87,498 \\
         & 94.70\% & 95.49\% & 90.68\% & 94.72\% & 95.12\% & 90.96\% \\
        Standby & 1,270 & 950 & 2,322 & 1,270 & 1,139 & 2,309 \\
        & 1.32\% & 0.99\% & 2.41\% & 1.32\% & 1.18\% & 2.40\% \\
        Off & 3,829 & 3,391 & 6,639 & 3,808 & 3,554 & 6,386 \\
          & 3.98\% & 3.53\% & 6.90\% & 3.96\% & 3.69\% & 6.64\% \\
        \midrule
        \multicolumn{7}{l}{\textit{Production Statistics}} \\
        H$_2$ Prod. (t) & 4,329.3 & 4,969.8 & 4,494.6 & 4,503.9 & 4,637.3 & 4,315.3 \\
        Elec. Consumed (GWh)   & 200.5   & 202.2   & 192.2   & 200.6   & 201.4   & 192.7   \\
        \midrule
        \multicolumn{7}{l}{\textit{Device Wear (count)}} \\
        Cold Starts        & 162 & 149 & 9     & 90    & 111 & 8     \\
        Warm Starts        & 1,270 & 950 & 2,322 & 1,270 & 1,139 & 2,309 \\
        Stack Replacements & 1   & 1   & 0     & 1     & 1   & 0     \\
        \bottomrule
    \end{tabularx}
\end{table}

The highest impact parameter is the initial electrical efficiency of the device. Intuitively, the higher hydrogen production throughput for these devices introduces more revenue opportunities at no cost. This is confirmed in Table \ref{tbl: tornado part1} where all operational statistics are relatively unchanged between the lower and upper bound besides the hydrogen production quantity. Approaching the DOE technical targets for AWE efficiency (74\% LHV, \citep{TechnicalTargetsLiquid2025}) is the most optimistic option for improving economic feasibility of this technology. 

The plot in Figure \ref{fig: tornado plot} displays the operating degradation of the electrolyzer was the second highest impact device parameter with a NPV range of \$0.77MM. This comes unsurprisingly with the high on operation fraction ($\sim$95\% from Table \ref{tbl: result summary}). With reduced operation degradation, the optimal device operation included nearly twice the amounts of off and standby operation (8,961hr) compared to the high degradation counterpart (5,078hr). This shift in operation schemes reflects a discrete decision to avoid replacement in the given 11 year time horizon. These insights are a unique product of the multi-scale framework presented. 

The electrolyzer stack cost sensitivity was observed with projections of cheaper electrolysis stack materials \citep{krishnanPresentFutureCost2023} and targets \citep{TechnicalTargetsLiquid2025}. We found a moderate impact to the economics of the system. At the lower stack cost of \$150/kW, replacement occurred at the same frequency with the higher NPV due to the discounted stack price. However, with the higher stack cost of \$350/kW, the optimal operation shifted to minimize overall degradation accumulation and avoid stack replacement. This is reflected in the decrease in total hydrogen production (4,315.3t vs 4,637.3t) and nearly double the amount of hours in 'standby'/'off' modes (4,693hr vs. 8,695hr combined). As with the operation degradation parameter, the higher stack cost effects through the discrete replacement decision rather than on the economics directly.

 \begin{table}[!htp]
      \centering
      \caption{Tornado sensitivity analysis results -- Part 2. The two least impactful parameters (smallest NPV range) displayed here.}
      \label{tbl: tornado part2}
      \small
      \begin{tabularx}{\textwidth}{X cc cc}
          \toprule
          \textbf{Metric} & \multicolumn{2}{c}{\textbf{Startup Degradation}} & \multicolumn{2}{c}{\textbf{Standby Load}} \\
          \cmidrule(r){2-3} \cmidrule(r){4-5}
          & \textit{0.0} & \textit{122.4} & \textit{1.0\%} & \textit{10.0\%} \\
          & \multicolumn{2}{c}{$\mu$V/start} & \multicolumn{2}{c}{\% Capacity} \\
          \midrule
          \multicolumn{5}{l}{\textit{Economic Performance}} \\
          NPV (\$MM) & 1.57 & 1.56 & 1.57 & 1.55 \\
          LCOH (\$/kg) & 1.99 & 1.99 & 1.99 & 1.99 \\
          \midrule
          \multicolumn{5}{l}{\textit{Operation Mode (hours, \% total)}} \\
          On & 91,377 & 91,480 & 91,395 & 91,562 \\
             & 94.99\% & 95.10\% & 95.01\% & 95.19\% \\
          Standby & 0 & 1,208 & 1,363 & 881 \\
                  & 0.00\% & 1.26\% & 1.42\% & 0.92\% \\
          Off & 4,816 & 3,505 & 3,435 & 3,750 \\
              & 5.01\% & 3.64\% & 3.57\% & 3.90\% \\
          \midrule
          \multicolumn{5}{l}{\textit{Production Statistics}} \\
          H$_2$ Prod. (t)      & 4,644.4 & 4,648.5 & 4,645.0 & 4,651.9 \\
          Elec. Consumed (GWh) & 201.0   & 201.4   & 201.1   & 201.6   \\
          \midrule
          \multicolumn{5}{l}{\textit{Device Wear (count)}} \\
          Cold Starts        & 1,416 & 32    & 11    & 265 \\
          Warm Starts        & 0     & 1,208 & 1,363 & 881 \\
          Stack Replacements & 1     & 1     & 1     & 1   \\
          \bottomrule
      \end{tabularx}
  \end{table}

A surprising result concluded that the startup degradation had a negligible effect on the net present value of the electrolysis device. We considered no startup degradation (0 $\mu$V/start) and three times the experimental startup degradation voltage degradation (122.4 $\mu$V/start). The resulting range of NPV was \$0.02M in 11 years. There was a similar trend with the standby load percentage (see Table \ref{tbl: tornado part2}) where the range of NPV was \$0.01M. The minor effect of the startup degradation and the standby load indicates that these device parameters provide a nearly negligible NPV contribution. This is reflected in the cold and warm start statistics across all test cases, where warm starts dominate (besides 0.0 $\mu$V/start startup degradation).

We performed a further case study by removing the standby operation mode to confirm the economic indifference of the standby mode. Table \ref{tbl: 2-state} shows the results of this two-state model and the effect on the economics of the electrolysis system. We find that the inclusion of the third standby operation mode is not a significant overall difference to the NPV (\$0.04M). The two-state trial resulted in a \$1.53M NPV over the same 11-year time horizon compared to \$1.56M. We find that the inclusion of the startup degradation consideration remains economically arbitrary in the long-term scope. From this, we conclude that DR operation schemes may aggressively follow DAM signals with little compromise to device durability. Despite this, the warm start from standby operation remains more frequent than the cold start across all trials and should be used predominantly.

\begin{table}[!htp]
    \centering
    \caption{Comparison of results from presented three-state and two-state electrolysis model for DAM participation in
    the hub average price profile.}
    \label{tbl: 2-state}
    \small
    \begin{tabularx}{\textwidth}{X c c}
        \toprule
        \textbf{Metric} & \textbf{3-State} & \textbf{2-State} \\
        & \textit{('on'/'off'/'standby')} & \textit{('on'/'off')} \\
        \midrule
        \multicolumn{3}{l}{\textit{Economic Performance}} \\
        NPV (\$MM) & \$1.56 & \$1.53 \\
        LCOH (\$/kg) & 1.99 & 2.00 \\
        \midrule
        \multicolumn{3}{l}{\textit{Operation Mode (hours, \% total)}} \\
        On Operation & 91,489 (95.11\%) & 92,010 (95.65\%) \\
        Off Operation & 3,614 (3.76\%) & 4,183 (4.35\%) \\
        Standby Operation & 1,090 (1.13\%) & 0 (0\%) \\
        \midrule
        \multicolumn{3}{l}{\textit{Production Statistics}} \\
        H$_2$ Production (t)  & 4,649.1 & 4,660.8 \\
        Electricity Consumption (GWh) & 201.4 & 202.4 \\
        \midrule
        \multicolumn{3}{l}{\textit{Device Wear (count)}} \\
        Cold Starts & 151 & 906 \\
        Warm Starts & 1,090 & 0 \\
        Stack Replacements & 1 & 1 \\
        \bottomrule
    \end{tabularx}
\end{table}

\subsubsection{Design Tradeoffs}

It remains difficult to solely increase individual device parameters when designing new electrolysis systems. Thus, we seek to construct the optimal tradeoff heatmap between the device parameters. This map can be used to identify electrolyzer design tolerances which maximizes the economic value. This acts as an approach to compensate the inability to selectively improve single device parameters and strike a balance between different device priorities. We run this for the same 11-year time horizon as before for the same reasons.

\begin{figure}[!htp]
    \centering
    \begin{subfigure}[b]{0.48\textwidth}
        \centering
        \includegraphics[width=\textwidth]{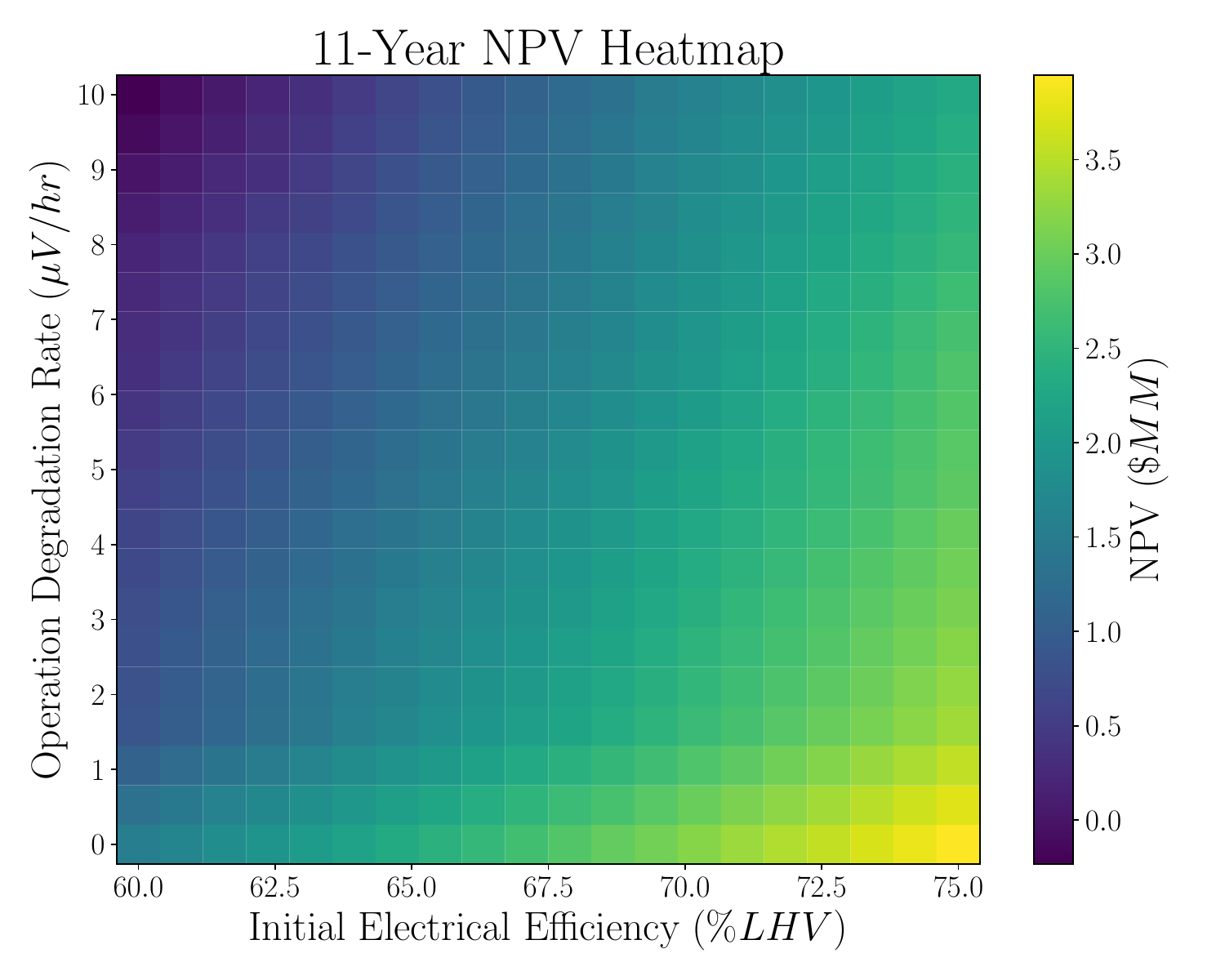}
        \label{fig:heatmap_npv}
    \end{subfigure}
    \hfill
    \begin{subfigure}[b]{0.48\textwidth}
        \centering
        \includegraphics[width=\textwidth]{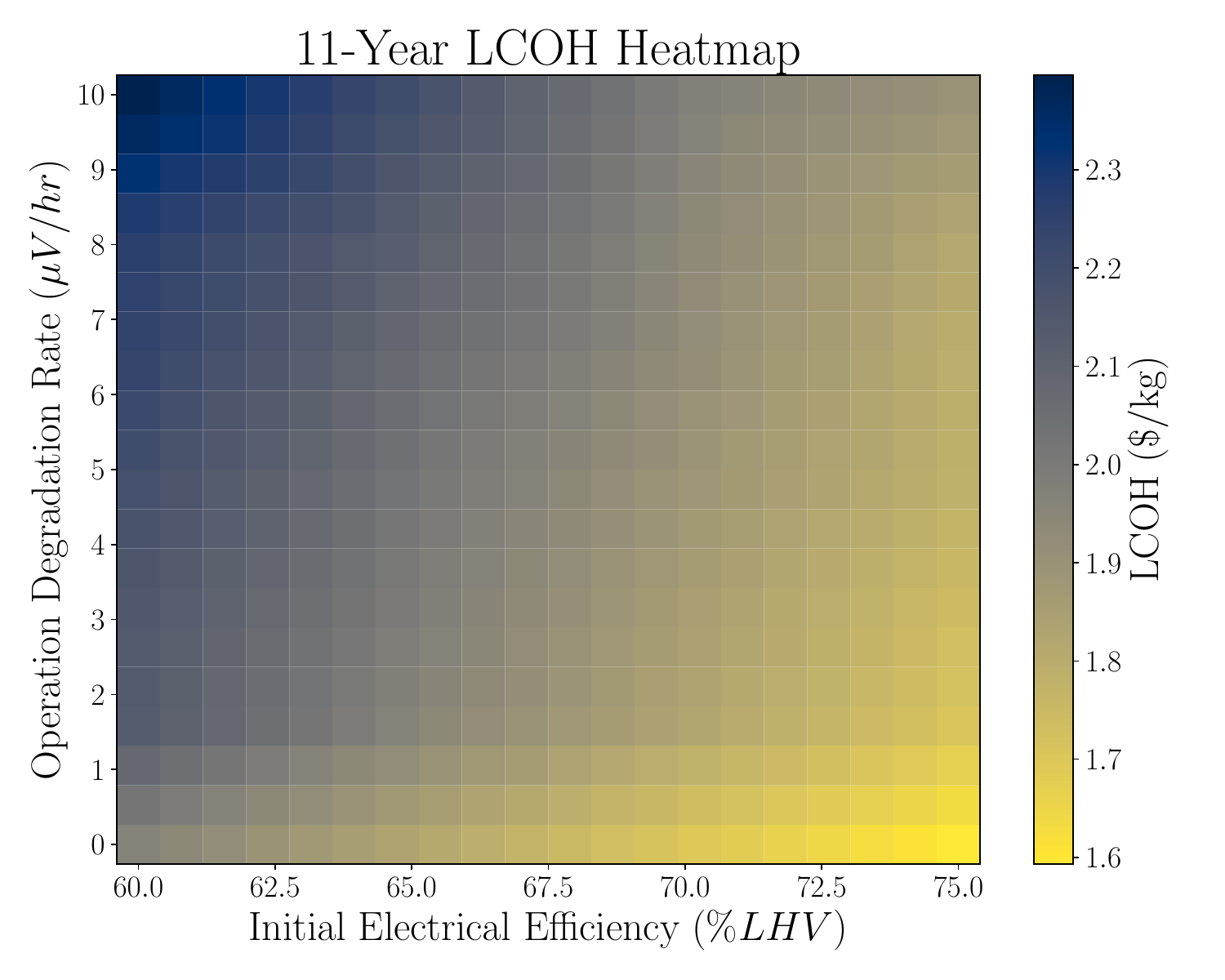}
        \label{fig:heatmap_lcoh} 
    \end{subfigure}
    \caption{Parameter heatmaps of 11-year NPV (left) and LCOH (right) as a function of initial electrical efficiency and operational degradation rate, using the ERCOT hub average LMP price profile (2014--2024). Values are linearly interpolated from 28 discrete grid points.}
	\label{fig: heatmap1}
\end{figure}

Figure \ref{fig: heatmap1} displays the parameter heatmap between initial efficiency and operation degradation for the 2.2 MW AWE considered previously. These plots demonstrate the NPV and LCOH of the electrolysis systems as a function of the different combinations of device parameters. This indicates where improvements in one parameter can compensate for deficiencies in the other for a set economic goal.

Both the NPV and LCOH are primarily sensitive to the initial electrical efficiency which coincides with the results from Figure \ref{fig: tornado plot}. Based on this, the DOE goal of 74\% LHV \citep{TechnicalTargetsLiquid2025} is the most direct methodology to increase economic viability. However, this also reveals that the LCOH fails to go below the overall goal of below \$1/kg even with the ideal no degradation and the DOE goal for electrical efficiency. In fact, the necessary CAPEX to drop below the general goal of \$1/kg LCOH given these electrolyzer design specifications is infeasible/negative. This reveals that without additional incentives and/or significant reductions in the cost of electricity: the goal of \$1/kg using electrolysis technology is infeasible within these device and market assumptions.

The sensitivity of NPV and LCOH to operation degradation is observed to be non-uniform across the design space. A linear relationship is observed between rates of 2-10 $\mu$V/hr. However, LCOH decreases from a sub-2 $\mu$V/hr degradation regime. For example, a 1 $\mu$V/hr operation degradation rate at 60\% LHV initial efficiency matches the LCOH for a 67.5\% LHV with an operation degradation rate of 9 $\mu$V/hr. This non-linear improvement in LCOH at low operation degradation rates is an manifestation of the transition in the optimal replacement schedule. The stack lifetime extends beyond the 11-year time horizon thus eliminating the significant replacement expenses.

\section{Conclusions and Future Work}\label{sec: conclusion}

We presented an optimization model of an electrolysis system to capture multi-scale interactions that link market participation, usage-based degradation, and stack replacement. We identified the long-term durability effects of flexibility provision for industrial scale electrolyzers and compared relative to static, non-flexible operation schemes. We also performed an economic sensitivity analysis to identify critical device design parameters and the optimal tradeoffs between devices. 

We conducted optimization over a representative case study of an AWE participating in three electricity price profile scenarios over a 22-year time horizon. First, we find that DR operation reduced electricity expenses by over 33\% for both ERCOT market participation scenarios compared to their respective static, non-flexible operation schemes. Across a 22-year horizon, the frequency of replacement remained constant across all cases of flexible and non-flexible operation. This cost decrease subsequently reduced the LCOH by over \$0.50/kg, approaching DOE production goals. Despite this, we found that given current ideal device specifications and market conditions, the goal of below \$1/kg LCOH is infeasible without supplementary incentives.

Investigating degradation effects, the rate of replacement for the flexible operation scheme demonstrated the capability to {\em extend} device electrolysis stack lifetime in highly volatile electricity price profiles.
We find that optimal DR operation during highly volatile markets hedged operation hours (load shifting) which consequently reduced the accumulation of relative degradation by up to 25\% in one year compared to static operation. 

To work towards implementing grid-integrated electrolysis, we performed a sensitivity analysis on device parameters and measured the effect on the NPV and LCOH of the system. Our findings indicate that minimizing steady-state operation degradation and maximizing hydrogen throughput are the most economically important parameters. In contrast, the values for the startup degradation and standby load percentage had a relatively negligible impact on the economic value of the electrolyzer when operating in DR. This demonstrates that optimal DR operation for electrolysis systems is not hindered by startup degradation accumulation. Rather, the steady-state operational degradation is one of the primary mechanism which influences profitable performance in electricity markets.

For further exploration of opportunities for industry to provide power grid flexibility, additional analyses are necessary to explore the impacts on durability of devices. A more practical model left for future work would include a rolling horizon market bidding formulation to better approximate active market participation. Additionally, future work to the proposed multi-scale optimization framework includes (1) participation in finer resolution electricity markets, such as the real-time market; (2) ancillary services provision, such as frequency regulation; (3) and implementing decomposition schemes to tractably solve larger instances of this hierarchical multi-scale optimization model. All three of these goals come from identifying additional pathways for providing flexibility in chemical manufacturing and identifying the effects on device durability. 

\section{Authorship Contribution}

\textbf{Kiernan X. Jennings: } Conceptualization, Methodology, Software, Investigation, Writing - original draft, Visualization
\textbf{Victor M. Zavala: } Conceptualization, Writing - review \& editing, Supervision, Funding acquisition.
\textbf{Styliani Avraamidou: } Conceptualization, Writing - review \& editing, Supervision, Funding acquisition.

\section{Declaration of Competing Interest}

We have no conflict of interest.

\section{Acknowledgments}

We acknowledge financial support from the U.S. National Science Foundation (CMMI-2328160).

\section{Declaration of Generative AI}

During the preparation of this work, the authors used Claude Code in order to produce figures presented. Data for the plots was produced independently of Generative AI. The authors reviewed and edited the content as needed and takes full responsibility for the content of the published article.

\bibliographystyle{elsarticle-harv}
\bibliography{references_1}
\newpage
\appendix

\section{Data Availability}

All data and code needed to reproduce the results are provided in \url{https://github.com/zavalab/JuliaBox/tree/master/MultiScaleElectrolyzer}.

\section{Nomenclature Table}

\begin{tabular}{@{}lp{0.75\textwidth}@{}}
	\multicolumn{2}{@{}l}{\textbf{Indices and sets}}\\[0.2em]
	\hline
	Symbol & Description \\
	\hline
	$m \in \mathcal{M}$ & Set of long-term (annual) time intervals $m$ \\
	$d \in \mathcal{D}$ & Set of medium-term (daily) time intervals $d$ \\
	$i \in \mathcal{I}$ & Set of short-term (hourly) time steps $i$ \\
	$t \in T$ & Lexicographic tuple set of $(m, d, i)$ for entire time horizon \\
    $t \in T^*$ & Subset of $T$ which includes the first time instant of each year $m$\\
	[1.0em]
\end{tabular}

\begin{tabular}{@{}lp{0.75\textwidth}@{}}
	\multicolumn{2}{@{}l}{\textbf{Parameters}}\\[0.2em]
	\hline
	Symbol & Description \\
	\hline
	$\psi_t^E$  & Electricity market price at time $t$ [\$/MWh]\\
	$\psi^H$ & Selling price of hydrogen [\$/kg] \\
	$\psi^{stack}$ & Cost of electrolyzer stack replacement [\$k/MW] \\
	$\psi^{CAPEX}$ & Fixed capital expense of electrolyzer construction [\$k] \\
	$\psi^{OPEX}$ & Fixed operational expense of electrolyzer operation [\$k/MW] \\
	$\phi$ & Nominal capacity of electrolyzer [MW] \\
	$\gamma$ & Length of time steps of electricity market [hr] \\
	$\eta^{op}$ & Operational voltage degradation [$\mu$V/hr] \\
	$\eta^{start}$ & Startup voltage degradation [$\mu$V/start] \\
	$\delta^{op}$ & Operational degradation efficiency loss [kg/(MWh $\cdot$ hr)] \\
	$\delta^{start}$ & Startup degradation efficiency loss [kg/MWh] \\
	$\alpha^0$ & Initial efficiency of the electrolyzer device  [kg/MWh] \\
    $\bar{\alpha}$ & Upper bound on the efficiency of the electrolyzer device  [kg/MWh] \\
	$\underline{\alpha}$ & Lower bound on the efficiency of the electrolyzer device [kg/MWh] \\
	$\beta$ & Intercept of hydrogen production curve [kg/hr] \\
	$\zeta^{sb}$ & Standby power consumption [\%] \\
	$V_0$ & Thermoneutral voltage of water electrolysis reaction [V] \\
	$V_i$ & Initial voltage of electrolysis cell [V] \\
	$LHV$ & Lower heating value of hydrogen [kWh/kg H$_2$]\\
	$HHV$ & Higher heating value of hydrogen [MWh/kg H$_2$] \\
	$\sigma$ & Chemical demand for operational subproblem [kg H$_2$] \\
	$\theta$ & Power ramping limit for electrolyzer [\% nominal capacity] \\
	[0.4em]
	
\end{tabular}

\begin{tabular}{@{}lp{0.65\textwidth}@{}}
	\multicolumn{2}{@{}l}{\textbf{Variables}}\\[0.2em]
	\hline
	Symbol & Description \\
	\hline
	Variables $\forall m \in \mathcal{M}$ & \\
	$z_m \in \{0,1\} $ & Binary variable defining replacement decision at year $m$ \\
	$v_{m} \in \mathbb{R}^+ $ & Auxiliary variable used in electrolyzer stack replacement decision linearization in year $m$\\
	[0.2em]
	Variables $\forall t=(m, d,i) \in T$ & \\
	$x_t^{on} \in \{0,1\}  $ & Binary variable defining 'on' operation of the electrolyzer at time $t$\\
	$x_t^{sb} \in \{0,1\} $ & Binary variable defining 'standby' operation of the electrolyzer at time $t$\\
	$x_t^{off} \in \{0,1\} $ & Binary variable defining 'off' operation of the electrolyzer at time $t$\\
	$x_t^{start} \in \{0,1\} $ & Binary variable defining if the electrolyzer undergoes cold startup at time $t$\\
	$h_t \in \mathbb{R}^+ $ & Hydrogen produced by electrolyzer at time $t$\\
	$e_t \in \mathbb{R}^+ $ & Electricity usage by electrolyzer at time $t$ \\
	$a_{t} \in \mathbb{R}^+ $ & Electrolyzer efficiency at time $t$\\
	$w_{t} \in \mathbb{R}^+ $ & Auxiliary variable used in hydrogen production linearization at time $t$\\
	[0.4em]
\end{tabular}

\section{Supplementary Case Study Parameters}

\begin{sidewaystable}[!htp]
	\centering
	\caption{Parameters for case study of 2.2 MW AWE participating in the DAM for 22-year time horizon. Table split between AWE device parameters and economic parameters.}
	\label{tbl: device params}
	\small
	\begin{tabular}{ c | c | c | c }
		\hline
		\textbf{Parameter} & \textbf{Symbol} & \textbf{Value}  & \textbf{Source}\\
		\hline
		Nominal capacity & $\phi$ & 2.2 MW & \cite{krishnanPresentFutureCost2023}\\
		Initial/maximum efficiency & $\alpha^0, \bar{\alpha}$ & 19.48 kg H$_2$/MWh & \cite{TechnicalTargetsLiquid2025}\\
        Intercept & $\beta$ & 9.66 kg H$_2$/hr & \cite{baumhofOptimizationHybridPower2023}\\
        Operational voltage degradation & $\eta^{op}$ & 3.2 $\mu$V/hr & \cite{TechnicalTargetsLiquid2025} \\
        Operational degradation & $\delta^{op}$ & 3.33E-5 kg H$_2$/(MWh $\cdot$ hr) & Derived \\
		Startup voltage degradation & $\eta^{start}$ & 40.8 $\mu$V/start & \cite{parkImpactDegradationEconomics2025}\\
        Startup degradation & $\delta^{start}$ & 4.25E-4 kg H$_2$/MWh& Derived \\
        Standby load & $\zeta^{sb}$ & 5\% & \cite{baumhofOptimizationHybridPower2023} \\
        Thermoneutral voltage & $V_0$ & 1.48 V & Known \\
		Initial cell voltage & $V_i$ & 1.9 V & \cite{parkImpactDegradationEconomics2025}\\
        \hline
        Hydrogen selling price & $\psi^H$ & \$3/kg & Assumed \\
        Hydrogen daily demand & $\sigma$ & 750 kg H$_2$/day & Assumed \\ 
        Discount factor & $\rho$ & 5\% & \cite{superchiDevelopmentReliableSimulation2023} \\
		Replacement cost & $\psi^{stack}$ & \$ 250 kUSD/MW &  \cite{krishnanPresentFutureCost2023}\\
		Total capital expense & $\psi^{CAPEX}$ &\$3,993 kUSD & \cite{aminahoTechnoeconomicAssessmentsElectrolyzers2025} \\
		Fixed operation expenses & $\psi^{OPEX}$ & 2\% CAPEX/yr & \cite{namiTechnoeconomicAnalysisCurrent2022} \\
		\hline
	\end{tabular}
\end{sidewaystable}

\section{Supplementary Sensitivity Analysis Results}

\end{document}